\magnification=\magstep1
\input amstex
\documentstyle{amsppt}

\NoBlackBoxes
\NoRunningHeads

\pagewidth{5.15truein}
\pageheight{8.00truein}



\topmatter
\title Hilbert series of quadratic algebras associated with
pseudo-roots of noncommutative polynomials
\endtitle
\author Israel Gelfand, Sergei Gelfand,Vladimir Retakh, Shirlei Serconek and Robert Lee Wilson
\endauthor
\address
\newline
I.~G., V.~R. and R.~W.: Department of Mathematics, Rutgers University, Piscataway,
NJ 08854-8019
\newline
S.~G.:
American Mathematical Society, P.O.Box 6248, Providence, RI 02940
{\rm and} Institute for Problems of Information Transmission, 19, Ermolova str.,
Moscow, 103051, Russia
\newline
S.~S.:
IME-UFG
CX Postal 131
Goiania - GO
CEP 74001-970 Brazil
\endaddress

\email
\newline
I.~G. : igelfand\@ math.rutgers.edu
\newline
S.~G. : sxg\@ams.org
\newline
V.~R. : vretakh\@ math.rutgers.edu
\newline
S.~S. : serconek\@math.rutgers.edu
\newline
R.~W. : rwilson\@math.rutgers.edu
\endemail

\endtopmatter
\centerline{\bf Abstract}
\medskip

The quadratic algebras $Q_n$ are associated with
pseudo-roots of noncommutative polynomials.  We compute 
the Hilbert series of the algebras $Q_n$ and of the dual algebras $Q_n^!$.
\bigskip 
\centerline{\bf Introduction}
\medskip

Let $P(x)=x^n-a_1x^{n-1}+\dots +(-1)^na_n$ be a polynomial over a
ring $R$.
Two classical problems concern the polynomial $P(x)$:
investigation of the solutions of the equation $P(x)=0$ and
the decomposition of $P(x)$ into a product of irreducible polynomials.

In the commutative case relations between these two problems are well
known:
when $R$ is  a commutative division algebra, $x$ is a central variable,
and the equation $P(x)=0$ has roots $x_1,\dots , x_n$, then
$$
P(x)=(x-x_n)\dots (x-x_2)(x-x_1). \tag 0.1
$$

In noncommutative case relations between the two problems are
highly non-trivial. They were investigated by Ore  \cite {O} and others.
(\cite {L} is a good source for references, see also the book \cite {GLR}
where matrix polynomials are considered.)
More recently, some of the present authors have obtained results \cite
{GR3, GR4, W}
which are important for the present work.  For a division algebra $R$,
I. Gelfand and V. Retakh
\cite {GR3, GR4} studied connections between the coefficients of $P(x)$
and a generic set of solutions $x_1,\dots , x_n$
of the equation $P(x)=0$. They showed that for any  ordering
$I=(i_1,\dots , i_n)$
of $\{1,\dots , n\}$ one can construct elements $y_k$ , $k=1,\dots , n$,
depending
on  $x_{i_1},\dots , x_{i_k}$ such that
$$
a_1=y_1+y_2+\dots +y_n,
$$
$$
a_2=\sum _{i<j} y_j y_i, \tag 0.2
$$
$$
\vdots
$$
$$
a_n = y_n\dots y_2y_1.
$$
\noindent
These formulas are equivalent to the decomposition
$$
P(t)=(t-y_n)\dots (t-y_2)(t-y_1) \tag 0.3
$$
where $t$ is a central variable. Formula (0.3) can be viewed as a
noncommutative
analog of formula (0.1). A decomposition of $P(x)$ for a noncommutative
variable
$x$ is more complicated (see \cite {GR4}).

The element $y_k$, which is defined to be the conjugate of $x_{i_k}$
by a Vandermonde quasideterminant involving $x_1,\dots ,x_k$, is a
rational
function in
$x_{i_1},\dots , x_{i_k}$; it is symmetric
in $x_{i_1},\dots , x_{i_{k-1}}$.
(Quasideterminants were introduced and studied in
\cite {GR1, GR2}.  We do not need the explicit formula for $y_k$ here.)
It was shown in \cite {W} that the polynomials in $y_k$ for a fixed ordering $I$ which
are symmetric
in $x_l$ can be written as polynomials in the
symmetric functions
$a_1,\dots , a_n$ given by formulas (0.2). 
Thus these are the natural noncommuttive symmetric functions.

It is convenient for our purposes to use the notation $y_k=x_{A_k,
i_k}$
where $A_k=\{i_1,\dots , i_{k-1}\}$ for $k=2,\dots , n$, $A_1=\emptyset $.
In the generic case there are $n!$ decompositions of type (0.3). Such
decompositions
are given by products of linear polynomials $t-x_{A,i}$ where
$A\subset \{1,\dots , n\}$, $i\in \{1,\dots ,n\}$, $i\notin A$. It is
natural
to call the elements $x_{A,i}$ {\it pseudo-roots} of the polynomial $P(x)$.
Note that elements $x_{\emptyset , i}=x_i$, $i=1,\dots , n$, are roots of
 the polynomial $P(x)$.

\smallskip
In \cite {GRW} I. Gelfand, V. Retakh and R. Wilson introduced
the algebra $Q_n$ of all pseudo-roots of a generic noncommutative
polynomial.
It is defined  by generators $x_{A,i}$,
$A\subset \{1,\dots , n\}$, $i\in \{1,\dots ,n\}$, $i\notin A$
and relations
$$
x_{A\cup i, j}+x_{A,i}-x_{A\cup j, i}-x_{A,j},\tag 04a
$$
$$
x_{A\cup i, j}\cdot x_{A,i}-
x_{A\cup j, i}\cdot x_{A,j}, \ \ i,j \in \{1,\dots , n\} \setminus
A.\tag 0.4b
$$
In \cite {GRW} a natural  homomorphism $e$ of  $Q_n$ into the free skew
field
$\Cal C$ generated by $x_1,\dots , x_n$ was constructed. We believe that
the map $e$ is an embedding.

We consider the algebra $Q_n$ as a universal algebra of pseudo-roots of
a noncommutative polynomial of degree $n$. Our philosophy is the
following: the algebraic operations of addition, subtraction and
multiplication are cheap,
but the operation of division is expensive.
For our problem we cannot use the ``cheap" free
associative algebra generated by $x_1,\dots , x_n$, but to use the
gigantic
free skew field $\Cal C$ is too expensive. So, we suggest to use an
``affordable intermediate" algebra $Q_n$.

\smallskip
Relations (0.4) show (see \cite {GRW}) that we may define a linearly
independent set
of generators
$$
r_A=x_{A \setminus \{a_1\},a_1}+x_{A \setminus \{a_1,a_2\},a_2}+ \dots +
x_{\emptyset,a_k}
$$
for all nonempty $A = \{a_1,...,a_k\}\subseteq \{1,\dots ,n\}$.
These generators satisfy the quadratic relations
$$
\{r(A)(r(A \setminus\{i\}) - r(A\setminus \{j\})) + (r(A \setminus\{i\}) -
r(A\setminus \{j\}))r(A \setminus \{i,j\})$$ 
$$- r(A\setminus \{i\})^2 + r(A \setminus \{j\})^2 | i,j \in A \subseteq
\{1,...,n\}\}.
$$

Another linearly independent set of generators in $Q_n$, 
$\{u_A| \emptyset \ne A \subseteq \{1,....,n\}\}$,
supersymmetric to $\{r_A| \emptyset \ne A \subseteq \{1,....,n\}\}$, was used in \cite {GGR} for a
construction of noncommutative algebras related to simplicial
complexes.

As a quadratic algebra  $Q_n$ has a dual quadratic algebra
$Q_n^!$,
see \cite {U}. A study of this algebra is of an independent interest.
In Section  5  we describe generators and relations for the algebra
$Q_n^!$.

\smallskip
In this paper we compute the Hilbert series of the quadratic
algebras
$Q_n$ and $Q_n^!$ . Recall that if $W = \sum_{i \ge 0} W_i$ is a 
graded vector space with $dim \ W_i$ finite for all $i$
then the Hilbert series of $W$ is
defined by
$$
H(W,t)=\sum_{i \ge 0}({\text dim}W_i)t^i.
$$

Any quadratic  algebra  $A$ has a natural graded structure
$A=\sum_{i \ge 0}A_i$ where $A_i$ is the span of all products of $i$
generators. If $A$ is finitely generated then  the subspaces $A_i$ are
finite-dimensional and the Hilbert series $H(A,t)$ of $A$ is defined.
Note that the Hilbert series $H(A^!,t)$ is also defined for the dual
algebra $A^!$.

Recall that if $A$ and $A^!$ are Koszul algebras  then
$H(A,t)H(A^!,-t)=1$ (see \cite {U}). The converse is not true but the counter-examples
are rather superficial (see\cite {P, R}).

The following two theorems,
which are the main results of this paper, show that the quadratic
algebras $Q_n$ satisfy this necessary condition for
the Koszulity of $Q_n$.

\proclaim {Theorem 1}
$$H(Q_n,t) = {{1-t} \over {1-t(2-t)^n}}$$

\endproclaim

\proclaim {Theorem 2}
$$H(Q_n^!,t) = {{1+t(2+t)^n} \over {1+t}}$$

\endproclaim

In the course of proving Theorem 1 we develop results (cf. Lemma 4.4)
which describe the
structure of $Q_n$ in terms of $Q_{n-1}$.  These results appear to be of
independent interest.
We use these results to compute
(Corollary 4.9) the Hilbert series of $Q_n$ in terms of the Hilbert series
of $Q_{n-1}$.
While proving Theorem 2 we determine (Proposition 6.4) a basis for the
dual algebra $Q_n^!$.

We begin, in Section 1, by recalling, from \cite{GRW}, the construction of
$Q_n$ (as a quotient of the
tensor algebra $T(V)$ for an appropriate vector space $V$) and developing
notation for certain important
elements of $T(V)$.  We also note that $Q_n$ has a natural filtration.  In
Section 2 we study the associated
graded algebra $gr \ Q_n$, obtaining a presentation for $gr \ Q_n$.  In
view of the basis theorem for $Q_n$
(in \cite{GRW}) it is easy to determine a basis for $gr \ Q_n$.  We next,
in Section 3, define certain
important subalgebras of $Q_n$ which we denote $Q_n(1)$ and $Q_n(\hat 1)$.
We show that the structures
of these algebras are closely related to the structure of $Q_{n-1}$.  In
Section 4 we use these facts
to prove Theorem 1 by induction on $n$. We then begin the study of the
dual algebra $Q_n^!$, recalling
generalities about the algebra and finding the space of defining relations
in Section 5 and
constructing a basis for $Q_n^!$ in Section 6.  The proof of Theorem 2,
contained in Section 7,
is then straightforward.

\bigskip
\head 1. Generalities about $Q_n$ \endhead

\medskip
The quadratic algebra $Q_n$ is defined in \cite{GRW}. 
Here we recall one presentation of $Q_n$ and develop some notation.   
Let $V$ denote the vector space over a field $F$ with basis 
$\{ v(A) \; | \;  \emptyset \neq A \subseteq \{ 1, \dots, n \} \}$ and  $T(V)$ denote the tensor algebra on $V$.  The symmetric group on $\{1,...,n\}$ acts on $V$ by $\sigma (v(A)) = v(\sigma(A))$ and hence also acts on $T(V)$.

Note that $$T(V) = \sum_{i \geq 0} T(V)_i$$
where 
$$T(V)_i = {\text{span}} \{ v(A_1) \dots v(A_i) \; | \;  \emptyset \neq A_1, \dots, A_i 
\subseteq \{ 1, \dots, n \} \}$$
 is a graded algebra.  Each $T(V)_i$ is finite-dimensional.  

Also, defining 
$$T(V)_{(j)} = {\text{span}} \{ v(A_1) \dots v(A_i) \; | \; i \geq 0,  
| A_1 | + \dots + | A_i | \leq j \}$$ 
gives  an increasing filtration 
$$F.1 = T(V)_{(0)} \subset T(V)_{(1)} \subset \dots    $$
of $T(V)$. 

Note that 
$$T(V)_{(j)} = \sum_{i \geq 0} T(V)_i \cap T(V)_{(j)} .$$ 

Let $ \emptyset \neq B \subseteq A \subseteq  \{ 1, \dots, n \}$ and write 
$B = \{ b_1, \dots, b_k \}$  where $b_1 > b_2 > \dots > b_k$.  Let $Sym(B)$ 
denote the group of all permutations of $B$.   
When convenient we will write $A \setminus b_1...\setminus b_k$ in place of 
$A \setminus \{b_1,...,b_k\}.$  Define $\Cal V(A:B)$ to be
 
$$  \sum_{\sigma \in Sym(B)}   
 \ sgn(\sigma) \; \sigma \{ v(A)v(A \setminus b_1)v(A \setminus b_1 \setminus b_2) 
\dots v(A \setminus b_1 \dots \setminus b_{k-1}) $$
$$ + \sum_{u=1}^{k-1} (-1)^u \{ v(A \setminus {b_1})...v(A \setminus {b_1} ...\setminus {b_{u-1}})
v(A \setminus {b_1} ... \setminus {b_u})^2v(A\setminus {b_1} ... \setminus {b_{u+1}}) ... $$ 
$$v(A \setminus {b_1} ... \setminus {b_{k-1}})\}$$ 
$$ + (-1)^k  v(A \setminus {b_1}) ... v(A \setminus {b_1}...\setminus {b_k})\}.$$

Let $\Cal Q = {\text{span}} \{ \Cal V(A:B) \; | \; B \subseteq A \subseteq  
\{1, \dots, n \}, |B| = 2 \}$ and let $<\Cal Q>$ denote the ideal in $T(V)$ 
generated by $\Cal Q$.  Denote the quotient $T(V)/<\Cal Q>$ by $Q_n$.  Since 
$\Cal Q \subseteq T(V)_2$, $Q_n$ is a quadratic algebra.  $Q_n$ is, of course, 
graded;
$$ Q_n = \sum_{i \geq 0} Q_{n,i}, $$ 
where 
$$  Q_{n,i} = (T(V)_i + <\Cal Q>) / <\Cal Q> . $$
Defining  
$$  Q_{n,(j)} = (T(V)_{(j)} + <\Cal Q>) / <\Cal Q>  $$
gives an increasing filtration 
$$F.1 = Q_{n,(0)} \subset Q_{n,(1)} \subset \dots    $$
of $Q_n$.
Note that 
$$   Q_{n,(j)} = \sum_{i \geq 0}   Q_{n,i} \cap  Q_{n,(j)}. $$

\noindent Let $r(A)$  denote $v(A) + <\Cal Q>$ and 
$ \Cal R(A:B)$ denote $ \Cal V(A:B) +  <\Cal Q>$.  

Note that if $|B| = 2$ then  $ \Cal R(A:B) = 0$ (in $Q_n$).

\bigskip
\head 2. The associated graded algebra $gr \ Q_n$ \endhead
\medskip
Let  $\Cal X = {\text{span}} \{  v(A)(v(A \setminus i) - v(A \setminus j))  \; | \; 
    i,j \in A  \subseteq  \{1, \dots, n \} \}$. 
$$  X_n = (T(V) + <\Cal X>) / <\Cal X> . $$

\noindent Let $x(A)$  denote $v(A) + <\Cal X>$.  

Note that  $X_n$ is graded 
$$ X_n = \sum_{i \geq 0} X_{n,i}, $$ 
where 
$$  X_{n,i} = (T(V)_i + <\Cal X>) / <\Cal X>  $$
and has an increasing filtration 

$$F.1 = X_{n,(0)} \subset X_{n,(1)} \subset \dots    $$
where 
$$  X_{n,(j)} = (T(V)_{(j)} + <\Cal X>) / <\Cal X> \; .  $$

A {\it string} is a finite sequence $\Cal B=(B_1,\dots, B_l)$
of nonempty subsets of $\{1,\dots ,n\}$. 
We call $l=l(\Cal B)$ the {\it length} of $\Cal B$ and 
$|\Cal B|=\sum _{i=1}^l |B_i|$ the {\it degree} of $\Cal B$.  Let 
$S$ denote the set of all strings. If ${\Cal B} = (B_1,...,B_l)$ 
and ${\Cal C} = (C_1,...,C_m) \in S$ define ${\Cal {BC}} = (B_1,...,B_l,C_1,...,C_m)$
and $x(\Cal B) = x(B_1)\dots x(B_l)$.
For any set $W \subseteq S$ of strings we will denote $\{x(\Cal B) \; | \; 
\Cal B \in W \}$ by $x(W)$.  Note that $S$ contains the empty string $\emptyset$.  Let $S^{\circ} = S \setminus \{ \emptyset \}$.  For any subset $U \subseteq S$ let $U^{\circ} = U \cap  S^{\circ}$.    

We recall from \cite{GRW}, the definition of $Y \subseteq S$.  Let   
$ \emptyset \neq A=\{a_1,\dots , a_l\} \subseteq  \{1,\dots ,n\} $ where $a_1>a_2>\dots >a_l$ and $j \leq |A|$. Then we write $(A:j) = 
(A, A \setminus a_1)  \dots  (A \setminus a_1 \setminus \dots \setminus a_{j-1})$, a string of length $j$.  

Consider the following condition on a string  $(A_1:j_1)\dots (A_s:j_s)  \in S$:  
$$\text{if} \;\; 2 \leq i \leq s \;\;  \text{and} \;\; A_i \subseteq A_{i-1}, \;\; {\text{then}} \;\;  
|A_i| \neq |A_{i-1}|-j_{i-1}. \tag 2.1 $$

Let $Y = \{ (A_1:j_1)\dots (A_s:j_s)  \in S \; | \; (2.1) \; \text{is satisfied} \}.$  
It is  proved in \cite{GRW} that $r(Y)$ is a basis for $Q_n$.   

Suppose $\Cal B=(B_1,\dots, B_l)$ is a string. Recall, from \cite{GRW}, 
that we may define by induction a sequence of integers
$n(\Cal B)=(n_1,n_2,\dots , n_t)$,
$1=n_1<n_2<\dots <n_t=l+1$, as follows:

$n_1=1$,

$n_{k+1}=\min (\{j>n_k|\ B_j\nsubseteq B_{n_k}\  {\text or}\ 
|B_j|\neq |B_{n_k}|+ n_k-j\}\cup \{l+1\})$, 

and $t$ is the smallest $i$ such that $n_i=l+1$. 

We call $n(\Cal B)$ the {\it skeleton} of $\Cal B$. 

Let  $\Cal B=(B_1,\dots, B_l)$ be  a string  with skeleton 
 $(n_1=1,n_2,\dots , n_t=l+1)$.

Define  ${\Cal B}^{\vee}$  to be the string  
$(B_{n_1},n_2 - n_1)(B_{n_2},n_3 - n_2) \dots   (B_{n_{t-1}},n_t - n_{t-1})$.

Note that  $l({\Cal B}^{\vee}) = l(\Cal B)$  and  
  $|{\Cal B}^{\vee}| = |\Cal B|$.  

\proclaim{Proposition 2.1 }  $ x(\Cal B)  =  x({\Cal B}^{\vee}) .$
\endproclaim

{\it Proof:} If $t = 1$ then $l = 0$ so $\Cal B = {\Cal B}^{\vee}$  is the 
empty string and $ x(\Cal B)  =  x({\Cal B}^{\vee}) = 1.$
Assume $t = 2$, so ${\Cal B}^{\vee} = (B_1,l)$.
We will proceed by induction on $l$. If $l = 1$ then 
$\Cal B = (B_1) = {\Cal B}^{\vee}$  so there is nothing to prove. 
If $l = 2$, then ${\Cal B} = B_1,B_1 \setminus i)$ for some $i$ 
and ${\Cal B}^{\vee} = (B_1, B_1 \setminus j)$ for some $j$.  Since 
$x(B_1)x(B_1 \setminus i) =  x(B_1)x(B_1 \setminus j)$ by the defining relations, the result holds in this case.  

Now assume $l > 2$ and that the result holds for all 
$\Cal C = (C_1,\dots, C_k)$  with skeleton $(1,k+1)$  and $k < l$.  
We have $\Cal B = (B_1,\dots, B_{l-l})(B_l)$ so 
$x(\Cal B) = x(B_1,\dots, B_{l-1}) x(B_l)$.  Since the skeleton of $(B_1,\dots, B_{l-1})$ is $(1,l)$  the induction assumption  applies and shows 
that  $x(B_1,\dots, B_{l-1}) = x(B_1,l-1)$.
Let $b$ denote the largest element of $B_1$.  Then since $(B_1,l-1) = 
(B_1)(B_1 \setminus  b ,l-2)$
we have $x(\Cal B) =  x(B_1,...,B_{l-1})x(B_l) = x(B_1,l-1)x(B_l)
=x(B_1)x(B_1 \setminus  b ,l-2)x(B_l)$.  
If $b \notin B_l$ the induction assumption shows that this is  
$x(B_1)x(B_1 \setminus  b ,l-1) = x(B_1,l)$  as required. 
So we may assume $b \in B_l$ and then,  since $|B_l| < |B_1|$, we may find 
$c \in B_1$, $c \neq b$, $c \notin B_l$.  Then by the induction assumption 
 $x(B_1,l-1)x(B_l) = x(B_1)x(B_1 \setminus  c ,l-2)x(B_l)$ and, again by the 
induction assumption, this is equal to  $x(B_1)x(B_1 \setminus  c ,l-1)$.  
 
Write $(B_1)(B_1 \setminus  c ,l-1) = (B_1, C_2, \dots, C_l)$ and note that, as \break $l > 2$, 
the largest element of $B_1$ is not in $C_l$.  Then by the previous 
case $x(B_1, C_2, \dots, C_l) = x({(B_1, C_2, \dots, C_l)}^{\vee})$.  
But $x(\Cal B) = x(B_1, C_2, \dots, C_l)$ and ${(B_1, C_2, \dots, 
C_l)}^{\vee} = (B_1,l)$ proving the result in case $t=2$.

Finallly, suppose $t > 2$ and suppose $n({\Cal B}) = (n_1,...,n_t)$.  
We proceed by induction on $t$.  \
Let ${\Cal B}' = (B_1,...,B_{n_2 - 1})$ 
and ${\Cal B}'' = (B_{n_2},...,B_l)$.  Note that $n({\Cal B}') = (1,n_2)$ 
and $n({\Cal B}'') = (n_2,...,n_t)$, and so, by induction, 
$x({\Cal B}') = x({\Cal B}'{^\vee})$ and $x({\Cal B}'') = x({\Cal B}''{^\vee})$.  
Then $x({\Cal B}) = x({\Cal B}')x({\Cal B}'') = x({\Cal B}'{^\vee})x({\Cal B}''{^\vee})
=x({\Cal B}^{\vee})$, proving the proposition.

\vskip 8 pt

Let $gr \; Q_n$ denote the associated graded algebra of $Q_n$.  For any string 
$\Cal B$  let $\bar r(\Cal B)$ denote the element $r(\Cal B) + 
Q_{n, |\Cal B|-1}$  of $gr \; Q_n$.  

For any set $S$ of strings write $\bar r(S) = \{ \bar r(\Cal B) \; | \; 
\Cal B \in S \}$.  

\proclaim{Lemma 2.2}  $\bar r(Y)$ is a basis for $gr \; Q_n$.

\endproclaim

{\it Proof:} This follows from the fact that  $r(Y)  \cap Q_{n,i}$ is a basis 
for  $Q_{n,i}$ (Theorem 1.3.8 and Proposition 1.4.1 of \cite{GRW}).  

\proclaim{Corollary 2.3}  The linear map  $\phi : X_n  \to  gr \; Q_n$ defined by 
$\phi(x(\Cal B)) =  \bar r(\Cal B)$  is an isomorphism of algebras.
\endproclaim

{\it Proof:} Since $Q_n$ is generated by $\{ r(A) \; | \; \emptyset \neq A 
\subseteq \{1, \dots, n \} \}$, $gr \; Q_n$ is generated by 
$\{ \bar r(A) \; | \; \emptyset  \neq A \subseteq \{1, \dots, n \} \}$.  
Since  $ 0 =  \Cal R(A: \{ i,j \}) = $ 
$$ r(A)(r(A \setminus i) - r(A \setminus j)) +  
(r(A \setminus i) - r(A \setminus j))r(A \setminus i \setminus j) - r(A \setminus i)^2 + r(A \setminus j)^2$$
we have $$ r(A)(r(A \setminus i) - r(A \setminus j)) \in Q_{n,2|A|-2}$$ 
and so $ \bar r(A)(\bar r(A \setminus i) - \bar r(A \setminus j))  = 0 $ in 
$gr \; Q_n$.  Consequently there is a homomorphism from $X_n$ into $gr \; Q_n$ that 
takes $x(A)$ into $\bar r(A)$.  Since the generating set 
$\{ \bar r(A) \; | \; \emptyset \neq A \subseteq \{1, \dots, n \} \}$ is 
contained in the image of this map, the map is onto.  Note that 
$Y = \{ \Cal B \; | \; \Cal B = {\Cal B}^{\vee} \}$.  Thus by Proposition 2.1, $X_n$ is 
spanned by $x(Y)$.  Since the image of this set is the linearly independent 
set $\bar r(Y)$, the map is injective.  

\medskip

\bigskip
\head 3. The subalgebras $Q_n(1)$ and $Q_n(\hat 1)$ \endhead
\medskip
Let $Q_n(\hat 1)$ denote the subalgebra of $Q_n$ generated by 
$\{ r(A) \; | \; \emptyset \neq A \subseteq \{2, \dots, n \} \}$.  Let 
$S(1) =  \{ \Cal B = (B_1, \dots, B_l) \in  S \; | \; 1 \in B_i 
\; \text{for all i} \}$, $S(1)^\dag = \{{\Cal B} = (B_1,...,B_l) \in S(1) | |B_i| > 1 \; \text{for all i} \}$  and 
$S(\hat 1) = \{ \Cal B = (B_1, \dots, B_l) \in  S \; | \ B_1, \dots, B_l 
\subseteq \{2, \dots, n \} \}$.  Let 
$Y(1) = Y \cap  S(1)$, $Y(1)^\dag = Y \cap S(1)^\dag$, and $Y(\hat 1) = Y \cap  S(\hat 1)$.  

Let $Y_{(n-1)}$ 
denote $\{(B_1,...,B_l) \in Y(1)| B_1,...,B_l  \subseteq \{1,...,n-1\}\}.$

\proclaim{Lemma 3.1}  $Q_{n-1}$  is isomorphic to $Q_n(\hat 1)$.  
\endproclaim 
{\it Proof:}  For any subset $A \subseteq \{1, \dots, n-1 \}$, let 
$A+1$ denote $  \{ a+1 \; | a \in A \}$, a subset of  $\{2, \dots, n \}$.  
Clearly there is a homomorphism from  $Q_{n-1}$  into $Q_n(\hat 1)$  that 
takes $r(A)$ into $r(A+1)$.  This map is injective since the ``$r(Y)$-basis" 
for $Q_{n-1}$ maps into a subset of $r(Y) \subseteq Q_n$.   Since the 
generators for  $Q_n(\hat 1)$ are contained in the image of this map, it is 
onto.

\proclaim{Corollary 3.2} $Y(\hat 1)$ is a basis for $Q_n(\hat 1)$.
\endproclaim

Let $Q_n(1)$ denote the subalgebra of $Q_n$ generated by 
$\{ r(A) \; | \; 1 \in  A \subseteq \{1, \dots, n \} \}$.  

\proclaim{Lemma 3.3} The map from $gr \; Q_n(\hat 1)$ into $gr \; Q_n(1)$  that takes 
$\bar r(A)$  into $\bar r(A \cup \{1\})$ is an injective homomorphism and 
$\bar r({Y(1)}^{\dag})$ is a basis for the image.
\endproclaim

{\it Proof:}  $gr \; Q_n(\hat 1)$ has generators 
$\{\bar r(A) \; | \; \emptyset \neq  A \subseteq \{2, \dots, n \} \}$ and relations 
$\{ \bar r(A)(\bar r(A \setminus i) - \bar r(A \setminus j))|i,j \in A \subseteq \{2,...,n\}\}$.  Since \break
$ \bar r(A \cup \{1\})(\bar r(A \setminus i \cup \{1\}) - 
\bar r(A \setminus j \cup \{1\})) = 0$ in  $gr \; Q_n(1)$  the required homomorphism 
exists.  Since the homomorphism maps  $\bar r(Y(\hat 1))$ injectively to 
 $\bar r({Y(1)}^{\dag})$, a subset of $\bar r(Y)$, the homomorphism is injective and   
 $\bar r({Y(1)}^{\dag})$ is a basis for the image.     

\proclaim{Lemma 3.4}
\vskip 4 pt
a) $\bar r(Y(1))$ is a basis for $gr \; Q_n(1)$.
\vskip 4 pt
b) $r(Y(1))$ is a basis for $Q_n(1)$.
\endproclaim

{\it Proof:} a) Since  $\bar r(Y(1)) \subseteq \bar r(Y)$  it is linearly 
independent.  Hence it is sufficient to show that $\bar r(Y(1))$ spans 
$gr \; Q_n(1)$.  But $gr \; Q_n(1)$ is spanned by the elements 
$\bar r(\Cal B)$ where  $\Cal B = (B_1,\dots, B_l)$, $ 1 \in B_1, \dots, B_l$.  By Proposition 2.1 
$\bar r(\Cal B)  =  \bar r({\Cal B}^{\vee})$  where  $(n_1, \dots , n_t)$  is 
the skeleton of $\Cal B$ and    
${\Cal B}^{\vee} =
(B_{n_1},n_2 - n_1)(B_{n_2},n_3 - n_2) \dots   (B_{n_{t-1}},n_t - n_{t-1})$.
Since $1 \in B_j$ for each $j$,  ${\Cal B}^{\vee} \in Y(1)$ giving the result. 

Part b) is immediate from a).

\bigskip

If $A$ and $B$ are algebras, let $A*B$ denote the free product of $A$ and $B$
(cf. \cite{B}, Chap. 3, \S 5, Ex. 6).  Thus there exist 
homomorphisms $\alpha: A \rightarrow A*B$ 
and $\beta: B \rightarrow A*B$ such that if $G$ 
is any associative algebra and $\mu: A \rightarrow G, \nu : B \rightarrow G$ 
are homomorphisms then there exists a unique homomorphism 
$\lambda: A*B \rightarrow G$ such that $\lambda \alpha = \mu$ and $\lambda \beta = \nu.$
Furthermore, if $A$ and $B$ have identity element $1$,  
$\{1\} \cup \Gamma_A$ is a basis for $A$ and $\{1\} \cup \Gamma_B$ is a basis for $B$
then $A*B$ has a basis consisting of $1$ and all 
products $g_1,...,g_n$ or $g_2,...,g_{n+1}$ where 
$n \ge 1$ and $g_t \in \alpha(\Gamma_A)$ if $t$ is even 
and $g_t \in \beta(\Gamma_B)$ if $t$ is odd.

\proclaim{Lemma 3.5} $gr \; Q_n(1)$  is isomorphic to $gr \; Q_{n-1}*F[\bar r(1)]$.
\endproclaim
{\it Proof:} 
Let $\alpha: gr \; Q_{n-1} \rightarrow gr \; Q_{n-1} * F[\bar r(1)]$ and 
$\beta: F[\bar r(1)] \rightarrow gr \; Q_{n-1}*F[\bar r(1)]$ be the 
homomorphisms occuring in the definition of $gr \; Q_{n-1}*F[\bar r(1)]$.

If $\emptyset \ne A = \{a_1,...,a_k\}\subseteq \{1,...,n-1\}$ 
define $$\delta(A) = \{1,1+a_1,...,1+a_k\}.$$  Then define a map
$\mu: \{\bar r(A) | \emptyset \ne A \subseteq \{1,...,n-1\}\} \rightarrow gr \; Q_n(1)$ by
$$\mu(\bar r(A)) = \bar r(\delta(A)).$$
In view of Lemma 2.2, $\mu$ extends to a linear map
$$\mu: gr \; Q_{n-1} \rightarrow gr \; Q_n(1).$$
By Corollary 2.3, $\mu$ preserves the defining relations
for $gr \; Q_{n-1}$ and so is a homomorphism.  Lemma 3.4 implies that $\mu$ is injective.
Note that $\bar r(1) \in gr \; Q_n(1)$ generates a subalgebra 
isomorphic to the polynomial algebra $F[\bar r(1)].$  Thus there is an injection
$$ \nu: F[\bar r(1)] \rightarrow gr \; Q_n(1).$$

Consequently there is a homomorphism 
$$\lambda: gr \; Q_{n-1} * \Gamma[\bar r(1)] \rightarrow gr \; Q_n(1)$$
such that $\lambda\alpha = \mu$ and $\lambda\beta = \nu$.  We claim
that $\lambda$ is an isomorphism.

Let $\Cal T$ denote the set of all strings ${\Cal G}_1 ... {\Cal G}_n$ 
or ${\Cal G}_2 ... {\Cal G_{n+1}}$ where
${\Cal G}_i = {\Cal B}_i \in Y(1)^\dag$ if $i$ is odd and ${\Cal G}_i = \{1\}^{j_i}$ if $i$ is even.
Note that ${\Cal T} \subseteq S(1).$  Define
$$\Phi: {\Cal T} \rightarrow Y(1)$$ 
by $\Phi({\Cal B}) = {\Cal B}^\vee$.  
Define $\Psi: Y(1) \rightarrow {\Cal T}$ by 
$\Psi ((A,j)) = (A,j)$ if $j < |A|, \Psi((A,j)) = (A,j-1)\{1\}$ if $j = |A|,$ 
and $\Psi ((A_1,j_1)...(A_s,j_s)) = \Psi((A_1,j_1))...\Psi((A_s,j_s))$ 
if $(A_1,j_1)...(A_s,j_s)$ satisfies (2.1).  Then $\Phi$ and $\Psi$ are inverse mappings.

Let $\gamma_i = \alpha$ if $i$ is odd and $\gamma_i= \beta$ if $i$ is even. 
Then $gr \; Q_{n-1}*F[\bar r(1)]$ has basis consisting of $1$
and all products $\gamma_1\bar r({\Cal H}_1)...\gamma_n\bar r({\Cal H}_n)$ or 
$\gamma_2\bar r({\Cal H}_2)...\gamma_n\bar r({\Cal H}_{n+1})$
where ${\Cal H}_i \in Y_{(n-1)}$ if $i$ is odd and ${\Cal H}_i = \{1\}^{j_i}$ if $i$ is even.
Then $\lambda(\gamma_1\bar r({\Cal H}_1)...\gamma_n\bar r({\Cal H}_n)) =
\bar r(\delta({\Cal H}_1))\bar r({\Cal H}_2)...) 
= \bar r(\delta({\Cal H}_1){\Cal H}_2 ...)
\break = \bar r((\delta({\Cal H}_1){\Cal H}_2 ...)^\vee)$ and 
$\delta({\Cal H}_1){\Cal H_2}... \in {\Cal T}.$
Also, $\lambda(\gamma_2\bar r({\Cal H}_2)...\gamma_{n+1}\bar r({\Cal H}_{n+1}) \break
= \bar r({\Cal H}_2\delta({\Cal H}_3) ...)
= \bar r({\Cal H}_2\delta({\Cal H}_3) ...)^\vee)$
and ${\Cal H}_2\delta({\Cal H}_3)... \in {\Cal T}.$
Every element of $\Cal T$ arises in this way.  
Since $\Phi: {\Cal T} \rightarrow Y(1)$ is a bijection, we see that $\lambda$ 
maps a basis of $gr \; Q_{n-1}*F[\bar r(1)]$ 
bijectively onto the basis $\bar r(Y(1))$ of $gr \; Q_n(1)$, proving the lemma.

\bigskip
\head4. Proof of Theorem 1 \endhead
\medskip
Let $\theta : S \times  S \to S$ be defined by 
$$\theta ((B_1, \dots, B_l),(C_1, \dots, C_k))  = (B_1, \dots, B_l, C_1, \dots, C_k).$$  

\proclaim{Lemma 4.1}  
If $\Cal B = (B_1, \dots, B_l)$, $\Cal C = (C_1, \dots, C_k) \in Y$, $1 
\notin B_l$ and $1 \in C_1$, then $\Cal B \Cal C \in Y.$     
\endproclaim

{\it Proof:} 
Since $\Cal B \in Y$ we may write $\Cal B (A_1,j_1) \dots (A_s,j_s)$ where 
(2.1) is satisfied.  Since $1 \notin B_l$ we have $1 \notin A_s$. Similarly since 
$\Cal C \in Y$ we may write $\Cal C = (D_1,m_1) \dots (D_t,m_t)$ where 
condition (2.1) holds. Since $1 \in C_1$ we have $1 \in D_1$.  Then   
$\Cal B \Cal C = (A_1,j_1) \dots (A_s,j_s)(D_1,m_1) \dots (D_t,m_t)$.   Since (2.1) holds 
for $\Cal B$ and $\Cal C$, and since 
$D_1 \nsubseteq A_s$ (for $1 \in D_1$, $1 \notin A_s$), (2.1) is satisfied 
for   $\Cal B \Cal C$ and so  $\Cal B \Cal C \in Y$.  

\vskip 4 pt
Let  $\Cal B = (B_1, \dots, B_l) \in S$.  Define 
$A(\Cal B)  = \{ i \; | \; 1 \leq i \leq l-1, 1 \in B_i, 1 \notin B_{i+1} \}$ 
and $a(\Cal B) = |A(\Cal B)|$.  Set $S_{\{ i \}} = \{ \Cal B  \in S \; | \; 
a(\Cal B) = i \}$.   Then $S$ is equal to the disjoint union 
$\bigcup_{i \geq 0} S_{\{ i \}}$.         

\vskip 4pt
\proclaim{Lemma 4.2} $(SS(1)^{\circ} \cap S_{\{ 0 \}}) \times (S(\hat 1)^{\circ}S 
\cap  S_{\{ i \}})$ injects into  $S_{\{ i+1 \}}.$ 
\endproclaim

{\it Proof:}
Let ${\Cal B} = (B_1,...,B_l) \in S_{\{i+1\}}.$  
There are $l+1$ pairs in $S \times S$ which $\theta$ maps to $\Cal B$, 
namely $(B_1,...,B_j) \times (B_{j+1},...,B_l)$ for $ 0 \le j \le l.$  Now $(B_1,...,B_j) \in S_{\{0\}}$ implies
$j \le min \ A({\Cal B})$ while $(B_1,...,B_j) \in SS(1)^\circ$ and $(B_{j+1},...,B_l) \in S_{\{i\}}$ implies
$j \in A({\Cal B})$ and the lemma follows.

Let $L(1) = S(1)^{\circ} \times  S(\hat 1)^{\circ}$, and  
$L(i+1)   = L(1) \times  L(i)$, $i \geq 1$.
   
\proclaim{Corollary 4.3} $\bigcup_{i \geq 0} S(\hat 1) \times L(i) \times  S(1)$ 
injects into $S$. 
\endproclaim

{\it Proof:} The $i^{th}$ term in the union maps into  $S_{\{ i \}}$, so it is enough to prove 
that this is an injection.  Write this term as 
$(S(\hat 1) \times  S(1)^{\circ}) \times  (S(\hat 1)^{\circ} \times  L(i-1) \times  S(1))$ 
and observe that the result follows by the lemma and by induction on $i$.  

\vskip 4 pt 
Let $M$ denote the span of $Y(1)^{\circ}Y(\hat 1)^{\circ} \cap Y$ and let $N$ denote the subalgebra of $Q_n$ generated by $M$.

\proclaim{Lemma 4.4}  
$N$ is isomorphic to the free algebra generated by $M$ and the map 
$$Q_n(\hat{1}) \otimes N \otimes Q_n(1) \rightarrow Q_n$$ induced by 
multiplication is an isomorphism of graded vector spaces.
\endproclaim

{\it Proof:}   Let  $W = Y(1)^{\circ}Y(\hat 1)^{\circ} \cap Y$ and let $W^i = 
W \times \dots \times W$ (i times).  Then $W$, being linearly independent, is a basis for $M$.  By Lemma 4.2, $\bigcup_{i \geq 0} W^i$
   injects into $S$.  Indeed, Lemma 4.1 shows that the image is in $Y$.   Thus     
$\bigcup_{i \geq 0} W^i$ injects onto a basis for $N$, so $N$ is isomorphic 
to the free algebra generated by $M$.  Again by Lemma 4.1 we have that 
$\bigcup_{i \geq 0} Y(\hat 1) \times  W^i \times  Y(1)$ maps into $Y$. Since any substring of an element of $Y$ is again in $Y$, this map is onto.  By Corollary 4.3  the map is 
an injection.  This proves the final statement of the lemma.  

\vskip 4 pt

We now recall some well-known facts about Hilbert series (cf. \cite{U}, \S 3.3).

\proclaim{Lemma 4.5} (a) If $W_1$ and $W_2$ are graded vector spaces 
then $$H(W_1 \otimes W_2,t) = H(W_1,t)H(W_2,t).$$

(b) If $W$ is a graded vector space , then $$H(T(W),t) = {{1} \over {1 - H(W,t)}} .$$

(c) If $A = \sum_{i \ge 0} A_i$ and $B = \sum_{i \ge 0}B_i$ are graded algebras
with $A_0 = B_0 = F.1$, then 
$${{1} \over {H(A*B,t)}} = {{1} \over {H(A,t)}} + {{1} \over {H(B,t)}} - 1.$$

\endproclaim

Let $U(A:j)$ denote the span of all strings  $(A_1:j_1)\dots (A_s:j_s)$ 
satisfying $(2.1)$ such that $1 \in A_i$ for all $i$ and  $(A_s:j_s) = (A,j)$. 

\proclaim{Lemma 4.6}
 
$$H(U(A:j),t) = t^j (1-t)^{n-|A|} H(Q_{n}(1),t).$$

\endproclaim    

{\it Proof:}  Since whenever the string
 $(A_1:j_1)\dots (A_{s-1}:j_{s-1})$  satisfies $(2.1)$ then the string
 $(A_1:j_1)\dots (A_{s-1}:j_{s-1})(\{1, \dots, n \},j)$  also satisfies $(2.1)$, we have  
$$H(U(\{1, \dots, n \},j),t) = t^j H(Q_{n}(1),t).$$

  We now proceed by downward induction on $|A|$, assuming the result is true whenever $|A| > l$.  Let $|A| = l$.  Then   
$$H(U(A:j),t) = t^j  H(Q_{n}(1),t) -  t^j \sum_{C \supseteq A, |C|=|A|+m, 
m \geq 1} H(U(C:m),t).$$  
By the induction assumption this is 
$$t^j ( 1 -   \sum_{C \supseteq A, |C|=|A|+m, 
m \geq 1} t^m (1-t)^{n-|C|} ) H(Q_{n}(1),t).$$

Let $C = D \cup A$ where $D \subseteq \{1, \dots, n \} \setminus A.$  Then the expression becomes 
$$t^j ( 1 -   \sum_{\emptyset \neq D \subseteq \{1, \dots, n \} \setminus A}  
 t^{|D|} (1-t)^{n-|A|-|D|} ) H(Q_{n}(1),t).$$
By the binomial theorem the quantity in parenthesis is $(1-t)^{n-|A|}$, 
proving the result.  

\bigskip

Let $B \subseteq  \{2, \dots, n \}$ and let $Z(B)$ denote the span of all strings
$(A_1:j_1)\dots (A_s:j_s)$ such that $1 \in A_1, \dots, A_s$, 
$(A_1:j_1)\dots (A_s:j_s)$ satisfies $(2.1)$,  $|B| = |A_s| - j_s$,  $A_s \supseteq B$.  

\proclaim{Lemma 4.7} 
$$ H(Z(B), t) =  t H(Q_{n}(1),t).$$
\endproclaim

{\it Proof:} Write $A_s = B \cup E \cup \{ 1 \}$ where $B \cap E = \emptyset$ 
and $E \subseteq  \{2, \dots, n \}.$   Then 
$$
Z(B) = \sum_{E \subseteq  \{2, \dots, n \} \setminus B}  U(B \cup E \cup \{ 1 \},|E|+1)
$$
and so 
$$
H(Z(B), t) =  \sum_{E \subseteq  \{2, \dots, n \} \setminus B} t^{|E|+1} 
(1-t)^{n-|B|-|E|-1} H(Q_{n}(1),t).
$$
By the binomial theorem this is $t H(Q_{n}(1),t).$

\vskip 8 pt

\proclaim {Lemma 4.8}

$$H(M,t) = (H(Q_n(\hat{1}),t) - 1)(H(Q_n(1),t) - 1)$$ 
$$ - tH(Q_n(1),t)(H(Q_n(\hat{1}),t) - 1).$$

\endproclaim

{\it Proof:} $M$ is the span of   $Y(1)^{\circ}Y(\hat 1)^{\circ} \cap Y.$  
The complement $\Cal Z$ of  $Y(1)^{\circ}Y(\hat 1)^{\circ} \cap Y$ in 
 $Y(1)^{\circ}Y(\hat 1)^{\circ}$ is the set of all strings 
$(A_1:j_1)\dots (A_s:j_s)(B_1, \dots , B_l)$ such that 
$(A_1:j_1)\dots (A_s:j_s) \in Y(1)^{\circ}$ satisfies $(2.1)$,
$|B_1| = |A_s| - j_s$,  $A_s \supseteq B_1$,
 $(B_1, \dots , B_l) \in Y(\hat 1).$  
 Let $Z$ denote the span of $\Cal Z$.   The lemma follows 
from  showing that  
$$ H(Z,t) = tH(Q_n(1),t)(H(Q_n(\hat{1}),t) - 1).$$  

\vskip 8 pt

For  $\emptyset \neq B \subseteq  \{2, \dots, n \}$  
let $P(B)$ denote the span of all strings in $\Cal Z$ 
such that $B_1 = B$ and $P_0(B)$ denote the span 
of all strings $(B_1, \dots , B_l) \in Y(\hat 1)$  such that $B_1 = B$.  Then 
$Z =  \sum_{\emptyset \neq B \subseteq  \{2, \dots, n \}} P(B)$ and 
$H(P(B),t) = H(Z(B),t)H(P_0(B),t).$  By Lemma 4.7, this is equal to 
$tH(Q_n(1),t)H(P_0(B),t)$.  Thus 
$$ H(Z,t) = \sum_{\emptyset \neq B \subseteq  \{2, \dots, n \}} H(P(B),t) $$
$$= \sum_{\emptyset \neq B \subseteq  \{2, \dots, n \}}  tH(Q_n(1),t)H(P_0(B),t)$$
$$= tH(Q_n(1),t) \sum_{\emptyset \neq B \subseteq  \{2, \dots, n \}} H(P_0(B),t).$$
But  $\sum_{\emptyset \neq B \subseteq  \{2, \dots, n \}} H(P_0(B),t) = H(Q_n(\hat{1}),t) - 1$  
and the lemma is proved.

\vskip 8 pt

\proclaim {Corollary 4.9}
$${1 \over H(Q_{n},t)} = (2-t)({1 \over H(Q_{n-1},t)}) -1.$$

\endproclaim

{\it Proof:}  $H(Q_{n},t) =   H(Q_{n}(\hat 1),t) H(N,t) H(Q_{n}(1),t)$ by Lemma 4.4.  
For brevity we write $H(Q_n(\hat 1),t) = a$ and $H(Q_n(1),t) = b.$
Then $$H(N,t) = {{1} \over {1 - H(M,t)}} =
{{1} \over {(1-t)b + a + (t-1)ab}}$$ 
and so
$${{1} \over {H(Q_n,t)}} =  {{(1-t)b + a + (t-1)ab} \over {ab}}
={{1-t} \over {a}} + {{1} \over {b}} + t-1.$$
Since
$${{1} \over {b}} = {{1} \over {a}} - t$$ this 
gives 
$${{1} \over {H(Q_n,t)}} = {{1-t} \over {a}} + {{1} \over {a}} - 1
={{2-t} \over {a}} - 1.$$  
Now $a = H(Q_n(\hat 1),t) = H(Q_{n-1},t)$ so the corollary follows.

Theorem 1 now follows from Corollary 4.9 and the fact that $Q_0 =  F.$

\bigskip

\head 5. Generalities about the dual algebra $Q_n^!$ \endhead
\medskip
Let $V^*$ denote the dual space of $V$.  Thus $V^*$ has basis $$\{v^*(A)| \emptyset \ne A \subseteq \{1,...,n\}\}$$ 
where $$<v(A),v^*(B)> = \delta_{A,B}.$$

Note that $T(V^*) = \sum_{i \ge 0} T(V^*)_i$ is a graded algebra where $T(V^*)_i 
= \ {\text {span}} \ \{v^*(A_1)...v^*(A_i)| \emptyset \ne A_1,...,A_i \subseteq \{1,...,n\}\}.$  
Also, $T(V^*)$ has a decreasing filtration 
$$T(V^*) = T(V^*)_{(0)} \supset T(V^*)_{(1)} \supset ... \supset T(V^*)_{(j)} \supset ...$$ where
$$T(V^*)_{(j)} = \ {\text {span}} \{v^*(A_1) ... v^*(A_i) | |A_1| + ... + |A_i| \ge j\}.$$
In fact $$T(V^*)_{(j+1)} = (T(V)_{(j)})^\perp$$ for $j \ge 0.$

Define $Q_n^! = T(V^*)/<{\Cal Q}^\perp>.$ 

We may explicitly describe ${\Cal Q}^\perp$ and thus give a presentation of $Q_n^!$. To do this 
define the following subsets of $T(V^*)_2$: 
$$S_1 = \{v^*(A)v^*(B) | B \not \subseteq A \ {\text {or}} \  |B| \ne |A|, |A| - 1 \},$$

$$S_2 = \{v^*(C)(\sum_{i \in C} v^*(C \setminus i)) + v^*(C)^2| |C| \ge 2\} ,$$

$$S_3 = \{(\sum_{i \notin C} v^*(C \cup i)v^*(C)) + v^*(C)^2| C \ne \{1,....,n\}\},$$

$$S_4 = \{s(\{1,...,n\})^2\}.$$

\proclaim {Theorem 5.1}  $S_1 \cup S_2 \cup S_3 \cup S_4$ spans $\Cal Q^\perp.$  
Therefore, $Q_n^!$ is presented by generators 
$\{v^*(A) | \emptyset \ne A \subseteq \{1,...,n\}\}$ 
and relations $S_1 \cup S_2 \cup S_3 \cup S_4.$
\endproclaim
\vskip 6 pt

Before beginning the proof of this theorem, 
we present some examples and develop some notation.
Set $$s(A) = (v^*(A) + <{\Cal Q}^\perp>)/<{\Cal Q}^\perp> \; \in \; Q_n^!.$$
Write $s(i)$ for $s(\{i\})$, $s(ij)$ for $s(\{i,j\})$, etc.

\bigskip
\noindent {\bf Example.}  (a) $Q_2^!$ is $5$-dimensional with basis
$$\{1,s(1),s(2),s(12),s(12)s(1)\}.$$

(b) $Q_3^!$ is $14$-dimensional with basis
$$\{1,s(1),s(2),s(3),s(12),s(13),s(23),s(123),s(123)s(12),s(123)s(13),$$ 
$$s(12)s(1),s(13)s(1),s(23)s(2),s(123)s(12)s(1)\}.$$
These assertions follow from Proposition 6.4.

\bigskip
Note that $Q_n^! = \sum_{i \ge 0} Q_{n,i}^!$ is graded where 
$$Q_{n,i}^! = (T(V^*)_i + <{\Cal Q}^\perp>)/<{\Cal Q}^\perp>$$
and that $Q_n^!$ has a decreasing filtration $$Q_n^! = Q_{n,(0)}^! \supseteq Q_{n,(1)}^! \supseteq ... \supseteq 
Q_{n,(j)}^! \supseteq ... $$ where $$Q_{n,(j)}^!= (T(V^*)_{(j)} + <{\Cal Q}^\perp>)/<{\Cal Q}^\perp>.$$  
Clearly $$Q_{n,(j)}^! = \sum_{i \ge 0} Q_{n,i}^! \cap Q_{n,(j)}^!$$ and
$$Q_{n,i}^! \cap Q_{n,(j)}^! = (0)$$ if $j > ni.$

Let 
$gr \ Q^!_n = 
\sum_{j=0}^{\infty} 
Q^!_{n,(j)}/Q^!_{n,(j+1)}$, 
the associated graded algebra of $Q^!_n$.  Denote $s(A) + Q^!_{n,(|A| + 1)} \in \ gr \ Q^!_n$ by $\bar s(A)$.  
Then $\{\bar s(A)|\emptyset \ne A \subseteq \{1,...,n\}\}$ generates $gr \ Q^!_n.$

{\it Proof of Theorem 5.1:}  We first show that each $S_h, 1 \le h \le 4$ is contained in $\Cal Q^\perp$, i.e., that $<{\Cal V}(A: \{c,d\}),u_h> = 0$ whenever $c < d, c,d \in A \subseteq \{1,...,n\}$ and $u_h \in S_h.$  For $h = 1$ or $4$ this is clear.  

If $h = 2$, we note that $$<{\Cal V}(A: \{c,d\}),v^*(C)(\sum_{i \in C} v^*(C \setminus i)) + v^*(C)^2> = 0$$ unless $A = C$, $A \setminus c = C$ or $A \setminus d = C$.  
In the first case,  
$$<{\Cal V}(A: \{c,d\}),v^*(C)(\sum_{i \in C} v^*(C \setminus i)) 
+ v^*(C)^2>$$ $$ = <v(A)v(A \setminus d) - v(A)v(A \setminus c),v^*(A)(\sum_{i \in A} v^*(A \setminus i))> = 0.$$  
In the second case $$<{\Cal V}(A: \{c,d\}),v^*(C)(\sum_{i \in C} v^*(C \setminus i)) 
+ v^*(C)^2> $$ $$= <-v(A \setminus c)v(A \setminus c \setminus d), v^*(A \setminus c)
(\sum_{i \in A \setminus c} v^*(A \setminus c \setminus i))>$$ $$ + <v(A \setminus c)^2,
v^*(A \setminus c)^2> = -1 + 1 = 0.$$  In the third case $$<{\Cal V}(A: \{c,d\}),v^*(C)
(\sum_{i \in C} v^*(C \setminus i)) + v^*(C)^2>$$ $$ = <v(A \setminus d)v(A \setminus c 
\setminus d), v^*(A \setminus d)(\sum_{i \in A \setminus d} v^*(A \setminus d \setminus i))> $$ 
$$+ <- v(A \setminus d)^2,v^*(A \setminus d)^2> = 1 - 1 = 0.$$

If $h = 3$, we note that $$<{\Cal V}(A: \{c,d\}),(\sum_{i \notin C} v^*(C \cup i)v^*(C)) + v^*(C)^2> = 0$$ 
unless $A \setminus c = C$, $A \setminus d = C$ or $A \setminus c \setminus d = C$.  In the first 
case $$<{\Cal V}(A: \{c,d\}),(\sum_{i \notin C} v^*(C \cup i)v^*(C)) + v^*(C)^2>$$ $$ 
= <-v(A)v(A \setminus c), \sum_{i \notin A \setminus c} v^*(A  \cup i \setminus c)v^*(A \setminus c)>$$ 
$$ + <v(A \setminus c)^2,v^*(A \setminus c)^2> = -1 + 1 = 0.$$  In the second case  
$$<{\Cal V}(A: \{c,d\}),(\sum_{i \notin C} v^*(C \cup i)v^*(C)) + v^*(C)^2>$$ 
$$ = <v(A)v(A \setminus d), \sum_{i \notin A \setminus d} v^*(A  \cup i \setminus d)v^*(A \setminus d)>$$ 
$$ + <-v(A \setminus d)^2,v^*(A \setminus d)^2> = 1 - 1 = 0.$$ In the third case 
$$<{\Cal V}(A: \{c,d\}),(\sum_{i \notin C} v^*(C \cup i)v^*(C)) + v^*(C)^2> $$ 
$$= <v(A \setminus d)v(A \setminus c \setminus d) - v(A \setminus c)v(A \setminus c \setminus d), $$
$$\sum_{i \notin A \setminus c \setminus d} 
v^*(A  \cup i \setminus c \setminus d)v^*(A \setminus c \setminus d)> = 0.$$  

We will now use downward induction on $l$ to show 
that $( S_1 \cup S_2 \cup S_3 \cup S_4) \cap T(V^*)_{(l)}$ 
spans $\Cal Q ^\perp \cap T(V^*)_{(l)}$ for all $l \ge 0$. 
Note that, since $T(V^*) = T(V^*)_{(0)}$, this will complete 
the proof of the lemma.  Now $\Cal Q^\perp$ is contained 
in $T(V^*)_2$ and $T(V^*)_2 \cap T(V^*)_{(2n+1)} = (0),$ 
so the result holds for $l = 2n+1$. Assume the result 
holds whenever $l > m$ and let $u \in {\Cal Q ^\perp} \cap T(V^*)_{(m)}.$  
Suppose that $m$ is even.  Then by subtracting an element in the span 
of $S_1$ we may assume that $ u \in \sum_{|C| = {m/2}} a_Cv^*(C)^2 + T(V^*)_{(m+1)}$ 
for some scalars $a_C$.  Then by subtracting an element in the 
span of $S_3 \cup S_4$ we may assume that $u \in T(V^*)_{(m+1)}$.  
Hence the induction assumption gives our result in this case.  
Now suppose that $m$ is odd. Then by subtracting an element in 
the span of $S_1$ we may assume that 
$u \in \sum_{|C| = {(m+1)/2}, i \in C} b_{C,i}v^*(C)v^*(C \setminus i) + T(V^*)_{(m+1)}$ 
for some scalars $b_{C,i}.$  Since $0 = <{\Cal V}(C:\{c,d\}),u>$ 
for all $C$ with $|C| = {(m+1)/2}$ we see that $b_{C,c} = b_{C,d}$ 
for all $c,d \in C$.  Then by subtracting an element in the 
span of $S_2$ we may assume that $u \in T(V^*)_{m+1}$.  
Hence the induction assumption gives our result in this 
case and the proof of the lemma is complete.

\bigskip
\head 6. A basis for $Q_n^!$ \endhead
\medskip
Let $B = \{b_1,...,b_k\} \subseteq A \subseteq \{1,...,n\}$ 
with $b_1 > ... > b_k$.  Define ${\Cal S}(A:B) \in Q^!_n$ by 
$${\Cal S}(A:B)  = s(A) s(A \setminus b_1)...s(A \setminus b_1 ... \setminus b_k).$$
Let $min \ A$ denote the smallest element of $A$.  
Define $${\Cal S} = \{{\Cal S}(A:B)| B \subseteq A \subseteq \{1,...,n\}, min \ A \notin B\},$$
$$\bar{\Cal S} (A:B) = {\Cal S}(A:B) + Q^!_{n,(1 + (|B|+1)(2|A|-|B|)/2)},$$ 
and $$\bar{\Cal S} = \{\bar{\Cal S}(A:B)| B \subseteq A \subseteq \{1,...,n\}, min \ A \notin B\}.$$ 
\proclaim {Lemma 6.1}   $\bar {\Cal  S} \cup \{\bar s(\emptyset)\}$ 
spans $gr \; Q^!_n$ and ${\Cal S} \cup \{s(\emptyset)\}$ spans $Q^!_n.$
\endproclaim 
\vskip 6 pt

{\it Proof}:  
It is sufficient to show the first assertion.  
We know that  $\{\bar s(A)|\emptyset \ne A \subseteq \{1,...,n\}\}$ 
generates $gr \ Q^!_n.$  Since the sets $S_1, S_3$ and $S_4$ are 
contained in ${\Cal Q}^\perp$ we see that $ \bar s(A) \bar  s(B) = 0$ 
unless $B \subset A$ and $|B| = |A| - 1$.  Furthermore, since the 
set $S_2$ is contained in ${\Cal Q}^\perp$, we 
have $\bar s(A)(\sum_{i \in A} \bar s(A \setminus i)) = 0$ 
for all $ A \subseteq \{1,...,n\}, |A| \ge 2.$  
Then if $i,j \in A \subseteq \{1,...,n\}$ we have 
$$\bar s(A) \bar s(A \setminus i) \bar s(A \setminus i \setminus j)$$ 
$$ = - \bar s(A)(\sum_{l \in A, l \ne i} \bar s(A \setminus l)) \bar s(A \setminus i \setminus j)$$ 
$$ = - \bar s(A) \bar s(A \setminus j) \bar s(A \setminus i \setminus j).$$  The lemma is then immediate.
\vskip 6 pt
\proclaim {Lemma 6.2}  Let $B = \{b_1,...,b_k\} \subseteq A \subseteq \{1,...,n\}$ with $b_1 > ... > b_k$ 
and $k > 2.$  Then, for $0 \le m \le k-2,$ 
$$<\sum_{\sigma \in \ Sym(B)} \ sgn(\sigma)\sigma \{ r(A \setminus {b_1}) ... 
r(A \setminus {b_1} ... \setminus {b_{k-1}})\},$$ 
$$
{V^*}^m{\Cal Q^\perp}{V^*}^{k-m-2}> = 0.$$
\endproclaim
\vskip 6 pt

{\it Proof:}  If $k = 3$ 
then $$ \sum_{\sigma \in \ Sym(B)} \ sgn(\sigma)\sigma \{ r(A \setminus {b_1}) ... r(A \setminus {b_1} ... \setminus {b_{k-1}})\}$$ 
$$= {\Cal V}(A:\{b_1,b_2\}) + {\Cal V}(A:\{b_2,b_3\}) + {\Cal V}(A:\{b_3,b_1\}) \in {\Cal Q},$$ 
so the result holds.  Now assume that $k > 3$ 
and that the result holds for $k-1.$  
Then it is sufficient to show that 
$$<\sum_{\sigma \in \ Sym(B)} \ sgn(\sigma)\sigma \{ r(A \setminus {b_1}) ... r(A \setminus {b_1} ... \setminus {b_{k-1}})\},$$ 
$$v^*(A \setminus b_1){V^*}^{m-1}{\Cal Q^\perp}{V^*}^{k-m-2}> = 0$$
whenever $m > 0$ and that 
$$<\sum_{\sigma \in \ Sym(B)} \ sgn(\sigma)\sigma \{ r(A \setminus {b_1}) ... r(A \setminus {b_1} ... \setminus {b_{k-1}})\},$$
$${\Cal Q^\perp}{V^*}^{k-3}v^*(A \setminus b_1 \setminus ... \setminus b_{k-1})> = 0.$$
Both of these are immediate from the induction assumption.
\vskip 6 pt
\proclaim {Lemma 6.3}  Let $B \subseteq A \subseteq \{1,...,n\}$ and $|B| = k \ge 2.$  
Then $${\Cal V}(A:B) \in \cap_{m = 0}^{k-2} V^mRV^{k-m-2}.$$
\endproclaim
\vskip 6 pt

{\it Proof:}  It is enough to show that $<{\Cal V}(A:B), {V^*}^m{\Cal Q^\perp}{V^*}^{k-m-2}> = 0$ 
for all $m, 0 \le m \le k-2.$ 
This is immediate if $k = 2$, since $R = \ {\text {span}} \ {\Cal Q}$.  
We will proceed by induction on $k$.  Thus we assume $k > 2$ 
and that the assertion is true for $k-1$.  
Now if $m > 0$ then ${V^*}^m{\Cal Q^\perp}{V^*}^{k-m-2} 
= \sum_{C \subseteq \{1,...,n\}} v^*(C){V^*}^{m-1}{\Cal Q^\perp}{V^*}^{k-m-2}.$
Now $v^*(C){V^*}^{m-1}{\Cal Q^\perp}{V^*}^{k-m-2}$ is orthogonal 
to ${\Cal V} (A:B)$ unless $A \subseteq C \subseteq A \setminus B$ 
and $|C| = |A|$ or $|A|-1$.  But if $|C| = |A| $
then $v^*(C){V^*}^{m-1}{\Cal Q^\perp}{V^*}^{k-m-2}$ 
is orthogonal to ${\Cal V} (A:B)$ by Lemma 6.2 and 
if $|C| = |A| - 1$,  the induction assumption 
yields the same result.  
Thus, the assertion holds if $m > 0$.  For $m = 0$ we must consider
${\Cal Q}^{\perp}V^{*k-2} = \sum_{C \subseteq \{1,...,n\}} {\Cal Q^\perp}{V^*}^{k-3}v^*(C)$.
Now, ${\Cal Q}^\perp{V^*}^{k-3}v^*(C)$ 
is orthogonal to   ${\Cal V} (A:B)$ unless $A \subseteq C \subseteq A \setminus B$ 
and $|C| = |A|-|B|$ or $|A|- |B| + 1$.  If $|C| = |A| - |B|$ then 
${\Cal Q}^\perp{V^*}^{k-3}v^*(C)$ is orthogonal 
to   ${\Cal V} (A:B)$ by Lemma 6.2 and if $|C| = |A| - |B| + 1$ the 
induction assumption yields the same result.  This completes the proof of the lemma.

\bigskip

Let $\Cal V$ denote the span of $\{{\Cal V}(A:B)| B \subseteq A \subseteq \{1,...,n\}\}$.  The lemma shows that $\Cal V$ is orthogonal to ${\Cal Q}^\perp$ and so, the pairing of $T(V)$ and $T(V^*)$ induces a pairing of $\Cal V$ and $Q^!_n$.

\proclaim {Proposition 6.4} ${\Cal S}$ is a basis for $Q^!_n.$
\endproclaim
\vskip 6 pt

{\it Proof:} 
Suppose $min \ A \notin B \subseteq A \subseteq \{1,...,n\}$.  Then 
$$<{\Cal V}(A:B \cup \{min \ A\}),{\Cal S}(A:B)> = 1$$ and 
$$<{\Cal V}(C,D),{\Cal S}(A:B)> = 0 $$ if $|C| < |A|$ or if $|C| = |A|$ 
and $C \ne A$ or if $min \ A \in D$ and $D \ne B$.  It is then easy to 
see that $\Cal S$ is linearly independent.  In view of Lemma 6.3, this proves the proposition.

\proclaim{Corollary 6.5}  $gr \; Q_n^!$ is presented by 
generators $\{v^*(A)|\emptyset \ne A \subseteq \{1,...,n\}\}$
and relations $\bar S_1 \cup  \bar S_2 \cup \bar S_3$ where
$$\bar S_1 = \{v^*(A)v^*(B) | B \not \subseteq A \ {\text {or}} \  |B| \ne |A|, |A| - 1 \},$$

$$\bar S_2 = \{v^*(C)\sum_{i \in C} v^*(C \setminus i)| |C| \ge 2\} ,$$

$$\bar S_3 = \{v^*(C)^2| \emptyset \ne C \subseteq \{1,....,n\}\}.$$
\endproclaim

\bigskip
\head 7. Proof of Theorem 2 \endhead
\medskip
Clearly if $i > 0$ then ${\Cal S} \cap Q^!_{n,i}$ is a basis for $Q^!_{n,i}$.  Thus, for $i > 0$, $dim \ Q^!_{n,i}$ is equal to  $|{\Cal S} \cap Q^!_{n,i}|$.  This is the same as $|\{{\Cal S}(A:B) \in {\Cal S}|\ |B| = i-1\}|.$  Now   
$|\{{\Cal S}(A:B) \in {\Cal S}|\ |B| = i-1 \ 
{\text {and}} \ |A| = u\}| = {n \choose u}{u-1 \choose i-1}$.  
Thus $$H(Q^!_n,t) = 1 + \sum_{i > 0}(\sum_{u=i}^n {n \choose u}{u-1 \choose i-1})t^i$$
$$ = 1 + t\sum_{i > 0}(\sum_{u=i}^n {n \choose u}{u-1 \choose i-1})t^{i-1}$$
$$ = 1 + t\sum_{v=0}^{n-1}(\sum_{u=v+1}^n {n \choose u}{u-1 \choose v})t^v$$
$$ = 1 + t\sum_{u=1}^n\sum_{v=0}^{u-1} {n \choose u}{u-1 \choose v}t^v$$
$$= 1 + t\sum_{u=1}^n{n \choose u}\sum_{v=0}^{u-1}{u-1 \choose v}t^v$$ 
$$ = 1 + t\sum_{u=1}^n{n \choose u}(t+1)^{u-1}$$
$$ = 1 + {t \over {t+1}}\sum_{u=1}^n{n \choose u}(t+1)^u$$ 
$$ = 1 + {t \over {t+1}}((2+t)^n - 1) = {1 \over {t+1}}(t+1 + t(2+t)^n - t)$$ 
$$ = {1 \over {t+1}}(1 + t(2+t)^n).$$
This completes the proof of Theorem 2.

\enddocument